\documentclass[conference]{IEEEtran}
\IEEEoverridecommandlockouts
\usepackage{cite}
\usepackage{amsmath,amssymb,amsfonts}
\usepackage{algorithmic}
\usepackage{graphicx}
\usepackage{textcomp}
\def\BibTeX{{\rm B\kern-.05em{\sc i\kern-.025em b}\kern-.08em
    T\kern-.1667em\lower.7ex\hbox{E}\kern-.125emX}}
\newcommand{\ie}{\textit{i}.\textit{e}.,~}
\newcommand{\eg}{\textit{e}.\textit{g}.,~}
\usepackage[ruled, noend, linesnumbered, vlined]{algorithm2e}

\usepackage{booktabs}
\usepackage{makecell}
\usepackage{url}
\usepackage{enumitem}
\usepackage[dvipsnames]{xcolor}  
\usepackage{xargs}   
\usepackage[colorinlistoftodos,prependcaption,textsize=small]{todonotes}
\usepackage{seqsplit}
\usepackage{times} 
\usepackage{graphicx} 
\usepackage{ amsthm} 
\usepackage{times,amsmath,epsfig}
\usepackage{epstopdf}
\usepackage{pifont}
\usepackage{comment}
\usepackage{subfig}

\begin{document}

\title{A multi-objective optimization framework for online ridesharing systems*}

\author{\IEEEauthorblockN	{Hamed Javidi}
    \IEEEauthorblockA{Dept. of Electrical Engineering \\ and Computer Science \\
     Cleveland State University\\
    Cleveland, USA\\
    Email: h.javidi@vikes.csuohio.edu}
    \and
    \IEEEauthorblockN{Dan Simon}
    \IEEEauthorblockA{Dept. of Electrical Engineering \\ and Computer Science \\
     Cleveland State University\\
    Cleveland, USA\\
    Email: d.j.simon@csuohio.edu}
    \and
    \IEEEauthorblockN{Ling Zhu}
    \IEEEauthorblockA{Ford Motor Company \\
    Ann Arbor, USA\\
    Email: lzhu40@ford.com}
    \and
	\IEEEauthorblockN{Yan Wang}
    \IEEEauthorblockA{ Ford Motor Company\\
    Cleveland, USA\\
    Email: ywang21@ford.com}

\thanks{This document is the results of the research project funded by Ford Motor company.}
}

\maketitle

\begin{abstract}
The ultimate goal of ridesharing systems is to match travelers who do not have a vehicle with those travelers who want to share their vehicle. A good match can be found among those who have similar itineraries and time schedules. In this way each rider can be served without any delay and also each driver can earn as much as possible without having too much deviation from their original route. We propose an algorithm that leverages biogeography-based optimization to solve a multi-objective optimization problem for online ridesharing. It is necessary to solve the ridesharing problem as a multi-objective problem since there are some important objectives that must be considered simultaneously. We test our algorithm by evaluating performance on the Beijing ridesharing dataset. The simulation results indicate that BBO provides competitive performance relative to state-of-the-art ridesharing optimization algorithms. 
\end{abstract}


\section{Introduction}
In recent years ridesharing systems have become extremely popular. They are becoming popular because commuters prefer to find a more economical way to take their trips. Ridesharing provides an appropriate alternative for urban transportation due to the potential benefits; \eg decreased traffic congestion, diminished fuel consumption, and reduced greenhouse gas emissions. Ridesharing could lead to more efficient use of empty car seats and could significantly alleviate some of the societal concerns about issues like traffic congestion and air pollution. 

People who do not have a can be serviced by a vehicle that has a similar route and overlapping travel times, or by a vehicle that can reroute to match the passenger's destination. 
There are multiple types of ridesharing systems; \eg taxi sharing, car sharing, courier services, scooter sharing, bike sharing. These ridesharing and shuttle services reduce costs for the riders and provides income for drivers. 

There are two categories of ridesharing: static and dynamic. In static ridesharing, all driver and rider requests are available at the beginning of the time period and no new riders or drivers appear later. In dynamic ridesharing, however, they can appear after the beginning of the time period. In other words, in static ridesharing all drivers and all riders and all corresponding data, including departure and arrival locations, and pickup and drop off times, are known at the beginning of the time period. In dynamic ridesharing new riders and drivers continuously appear to request rides and to provide transportation.

Static ridesharing is not realistic and cannot be used in a real-world application, so we focus on dynamic ridesharing. Due to the unpredictable nature of the problem, future driver appearances and rider requests are unknown.

The ridesharing problem \cite{b1, b2} is related to the vehicle routing \cite{b3} and multi-vehicle pickup and delivery problems \cite{b4, b5}, in which customers request to be picked up from their origins and dropped off at their destinations while satisfying vehicle capacity and time constraints. Effective and efficient optimization technology that matches drivers and riders is one of the necessary components for an optimized ridesharing system. There are some significant challenges in this problem, such as fast computation, scalability, seat utilization, and quality of service (e.g., average trip delay).

There are multiple objectives in this problem. One objective is minimizing the distance traveled by each rider and each driver. Another objective is minimizing the waiting time of the riders before being picked up by a driver. Another objective is maximizing the total benefit obtained through successful matches, which is defined based on the drivers' preferred mode of payment \cite{b6}. Another objective is matching rate, i.e., maximizing the number of riders who can be picked up and delivered to their destination. In this paper, we mainly focus on two important objectives that are well-known among researchers: matching rate and total travel distance.

There have been various approaches to solving the dynamic ridesharing problem. In \cite{b7}, a rolling horizon solution is proposed. This approach  periodically matches unmatched riders with drivers who have empty seats. At each iteration of the rolling horizon, a matching problem is solved with an objective function aimed at minimizing the total travel distance and maximizing the matching rate. There are several approaches to determine the frequency of the iterations in a rolling horizon framework.The main two approaches include periodic optimization with a fixed time step, and event-driven optimization, where an event can include the appearance of one or more new riders. Event-driven approaches are usually used when systems are expected to respond quickly to an event in the environment while periodic optimization may imply longer response delay \cite{b8}. 

The ridesharing problem is known to be a non-deterministic polynomial time (NP-hard) problem \cite{b9, b10}. That is, to find the optimal solution of the problem may require an enormous amount of computation time, depending on the problem size. Because of the nature of this problem, optimal solutions are unattainable in large-scale problems in a reasonable amount of time. There is a trade-off in this type of problem: optimal solution, or fast computation with near-optimal solution. A near-optimal solution can be found with an online algorithm whose computation time is fairly low. Research has been done on both approaches to this problem. In some cases, it is impractical to generate optimal solutions within a given (usually small) time bound \cite{b10, b11}. Furthermore, for some cases, achieving a good feasible solution in a short time is more desirable than finding the optimal solution. Consequently, a trade-off between optimality and computation time must be considered.

Early work on this problem focused on traditional integer programming approaches, which are limited to small scale problems; \eg 8 drivers and 96 riders \cite{b12, b13, b14}. This method returns the optimal solution for the problem by searching all possible states of the problem. However, the computation time of such methods is too high so this method can be applied only to small problems.

Heuristic approaches \cite{b11, b7, b21} have been proposed to solve the real-time taxi ridesharing problem. One approach is to find a greedy solution, which assigns one rider request at a time to the best available vehicle. Although these approaches are fast, the solution is not good enough. 
One of the most well-known algorithms that uses an approximation method to approach this problem is T-share \cite{b11}. They partitioned the Beijing road network using a grid and define a road network node as an anchor point of each cell and then used precomputed travel distances and travel times of the shortest paths between each pair of cells. The strategy is greedy so it tries to find the best vehicle match in terms of minimized travel distance for every rider's request. Since this work uses approximation in terms of time and distance it is not well-optimized. However, despite its drawbacks, since it does not need to find the shortest path for every rider from an entire map graph, it is very fast and can be used for online ridesharing.
 
In this paper, we propose an algorithm that leverages an evolutionary approach to search for the optimal solution to the ridesharing problem. Our contributions are as follows.
\begin{enumerate}[label=(\Alph*)]
\item \textbf{Developing a simulator.} We started our work using an existing simulator called Cargo \cite{b15}. This simulator is designed for the taxi ridesharing problem. In order to have a peer-to-peer ridesharing simulator we modified it from taxi-sharing to peer-to-peer ridesharing.
\item \textbf{Proposing an efficient matching algorithm.} We propose a quick and efficient matching algorithm based on an evolutionary algorithm that can be used for dynamic ridesharing.
\item \textbf{Implementing state-of-the-art papers.} In order to compare the performance of our algorithm with other research, we implemented state-of-the-art optimization methods for the dynamic ridesharing problem.
\end{enumerate}
\section{Problem statement and formulation}
In this section, we first define terms, including parameters and variables, then we present the mathematical formulation of the proposed ride-sharing system and introduce the objective functions.

\subsection{Definitions}

\textbf{Definition 1 (Route)}
Let $G = (V, E)$ be a a graph of a road network system, where $V$ is the set of vertices representing intersections and $E$ is the set of edges representing streets. Each edge $ (u, v) \in E $ has a weight $ w_{uv} $ indicating travel distance along the street. A route is a set of connected edges, and its cost is the sum of the weights of all the edges in the route.

\textbf{Definition 2 (Rider)}
A rider wants to be matched with a driver who can meet the rider's constraints. Each rider $ r $ begins at origin $ r_o $ and requests to be transported to destination $ r_d $ by a driver. Each rider is associated with a time window $ (r_e, r_l) $. The early time $ r_e $ is the earliest possible time the rider $ r $ can be picked up, and the late time $ r_l $ is the latest time that the rider can be dropped off. Riders requests will appear in the ride-sharing system in real time. 

\textbf{Definition 3 (Driver)}
A driver can be assigned to riders whose constraints match the driver's constraints. $D$ is a set of drivers. Each driver $ d \in D $ begins its trip at its origin $ d_o $ and ends at its destination $ d_d $. Each driver has its own time constraint $ (d_e,d_l) : d_e \le d_l $, where $ d_e $ is the earliest possible departure time from $ d_o $ and $ d_l $ is the latest possible arrival time at $ d_d $. Each $d \in D$ has properties $ d_{loc}$, $d_{cap}$, $d_{load}$, $d_{speed}$, $d_{sch} $. Here, $ d_{loc} \in V $ is the driver's current location, which is the current vertex of  driver $ d $; $ d_{cap} $ is the maximum seat capacity; $ d_{load} $ is the number of riders; $ d_{speed} $ is the driver's speed, which is assumed to be constant; and $ d_{sch} $ is the trip schedule of the driver $d$, which is defined below.

\textbf{Definition 4 (Query)} A query is a rider’s request for a ride and is a tuple $(r_o, r_d, r_e, r_l) $, each of which are defined above. $ r_r $ is the query's timestamp indicating when the query was submitted. It is assumed that the request time $r_r$ is equal to $r_e$.

\textbf{Definition 5 (Schedule)} Each driver has its own schedule. A schedule is a temporally ordered sequence of pickup and delivery points corresponding to $ n $ queries. Note that schedules dynamically
change over time as riders appear. A driver can have any number of riders up to its capacity. Each schedule $ d_{sch} = (d_o, p_1, p_2, . . . , p_n, d_d) $ is a sequence of points, where $ d_o $ is the driver’s origin, and $ p_k$, $(1 \leq k \leq n) $, is an origin or destination of a rider. $ d_d $ is the driver's destination. If a driver is assigned to a rider, the driver will move to $ r_o $ to pick up the rider, then to $ r_d $ to drop off the rider. The first item in each schedule is the driver's origin, the last item is the driver's destination, and each rider's origin and destination appear in the schedule such that $ r_o $ precedes the corresponding $ r_d $.

\textbf{Definition 6 (Minimal shortest path, MSP)}
Minimal shortest path is a function that returns a set of edges that defines the shortest driving distance between a given pair of points. This definition can be extended to a schedule with more than two points by concatenating the pair-wise MSPs in the order given in the schedule.

\textbf{Definition 7 (Distance overhead)} Each driver $ d \in D $ has an origin $ d_o $ and destination $ d_d $ so its minimal shortest path (MSP) is defined based on the shortest path from $ d_o $ to $ d_d $ and its minimal schedule will be $ d_{sch} = (d_o, d_d) $. If $r \in R$ is served by $d$, the trip schedule must be changed such that  $ d_{sch'} = (d_o, r_o, r_d, d_d)$. The distance overhead $ d_{dov} = \mbox{dist}(d_{sch'}) - \mbox{dist}(d_{sch}) $ is defined as the difference between the travel distance of the updated trip schedule and the driver's MSP.

\subsection{Objectives}

Various objectives have been used in the literature for the ridesharing problem. Ridesharing solutions usually rely on single-objective optimization. In this work we aim to use a multi-objective cost function. The computational time required to update drivers' schedules when rider requests appear should be considered when evaluating a solution method, but this metric is unfortunately excluded by most previous research \cite{b16}. The performance metrics that have been used most widely include the following.

\textbf{Matching rate.} The matching rate $ M_R $ is defined as the total number of riders in the matched rider set $ R_s $ (i.e., the number of riders who were matched with a driver to get a ride) divided by the total number of riders $ |R| $ who requested a ride.

\begin{equation} \label{eq:M_R}
M_R = \frac {|R_s| }{|R|}
\end{equation} 

\textbf{Average distance overhead.} This is the average distance overhead of all drivers.

\begin{equation} \label{eq:d_{ov}}
D_{ov} = \sum_{i=1}^{|D|} d_{dov}(i)
\end{equation}

\textbf{Average matching delay.} This is the average time between when a rider sends its ride request and when it gets matched. This factor quantifies how fast the solution method can find a match for a request.

\textbf{Multi-objective cost function.} As noted above, the ridesharing problem can include various objectives. We define a multi-objective function to solve this problem. We want to define an objective function that satisfies the greatest possible number of ride requests while minimizing the total distance overhead. We define cost as a weighted combination of $M_R$ and total distance overhead $D_{ov}$. We define the cost so that our problem is a minimization problem. $\alpha$ is a constant to define the weight of each part of the objective function.
\begin{equation} \label{eq:cost}
\mbox{Cost} =  \frac{\alpha  \sum_{i=1}^{|D|} d_{dov}(i)}{ \sum_{i=1}^{|R|} \mbox{MSP}(r_o(i),r_d(i))} + (1-\alpha)(1 - M_R) 
\end{equation} 

\subsection{Problem Formulation}

We define the ridesharing problem as follows.

\textbf{Definition 8 (Dynamic Ridesharing)} Given a set of drivers $D$ and a set of riders $R$ on a road network, the dynamic ridesharing problem is the task of finding the optimal driver-rider combinations such that the user-defined combination of Equation~\ref{eq:cost} is minimized, subject to the following constraints.
\begin{enumerate}[label=(\Alph*)]
\item \textbf{Capacity constraint.} The number of riders served by any driver $d$ cannot exceed the corresponding maximum seat capacity:

\begin{equation} \label{eq:cap_constraint}
d_{load} \leq d_{cap} \qquad \forall \;  d \in D
\end{equation} 

\item \textbf{Waiting time constraint.} The time for the driver to pick up the rider after receiving the request cannot be greater than the rider's waiting time constraint:

\begin{equation} \label{eq:wating_constraint}
r_l - \mbox{MSP}(r_o,r_d) r_{speed} \geq 0 \qquad \forall \;  r \in R
\end{equation} 

\item \textbf{Driver time constraint.} Each driver has own destination $ d_d $ as well as its own maximum arrival time $ d_l $, and this constraint must be satisfied when a rider is added to the driver's schedule.

\end{enumerate}

\section{Ridesharing Approaches}
Proposed approaches for matching and routing in ridesharing problems are categorized in two classes: local optimization and global optimization.

\subsection{Local Optimization}

Providing ridesharing at a city scale is a complex problem. One major challenge is to develop a real-time matching algorithm that can quickly determine the best driver to satisfy incoming ride requests. Local optimization approaches use a local search algorithm in a first-come first-serve way to greatly reduce the complexity of ridesharing optimization \cite{b11, b17, b18, b19}. The first request (the first $r \in R$) will be served by searching potential drivers that are feasible candidates for giving a ride to the rider. Then, using a greedy approach, the driver candidate who has minimum distance overhead $d_{dov}$ will be selected from among other candidates. These approaches usually provide a trip as quickly as possible but riders would not mind having some delay if the ride offer would be cheaper.

Greedy approaches \cite{b22, b10, b20} find a driver $d \in D$ who has a minimum $d_{dov}$ for each rider request. In other words, they calculate the overhead distance between the driver's current schedule and the new schedule, and put this cost in a priority queue. The second step of these algorithms is the validation step, which checks that the constraints for the rider, the driver, and riders who are already assigned to the driver are satisfied. Greedy approaches use a first-come first-serve (FCFS) method that only focuses on local optimization instead of reducing the overall travel cost. 

T-Share \cite{b11} is a greedy algorithm that focuses on building an efficient service system to quickly handle the incoming flow of queries using local optimization. It reduces the problem size and continuously re-optimizes the schedules. However, it only focuses on local optimization instead of reducing the overall travel cost.

Recently, Huang et al. \cite{b18} proposed a greedy algorithm that uses a leveraging kinetic-tree structure that is designed for efficient scheduling of dynamic requests and that adjusts routes on the fly.

The nearest neighbor (NN) algorithm \cite{b22} is a heuristic greedy algorithm that uses Euclidean distance to match potential drivers with riders. This algorithm is implemented using a priority queue to rank drivers base on their distance from the rider's origin. 
This algorithm consists of two steps. In the first step, the nearest driver is selected based on the distance between its current location and the rider's location. In the second step, validation, the rider is inserted in the selected driver's schedule if all constraints for the rider, the driver, and riders who are already assigned to this driver are satisfied. The first valid driver from the queue will be assigned to the rider.
The complexity of the NN algorithm is $O (D \log D)$ for $D$ drivers.

\subsection{Global Optimization}
Approaches in this category are mainly based on two strategies: approximation and heuristics. One algorithm, simulated annealing, has already been published, and another algorithm, BBO, is discussed in Section \ref{Sec:Proposed method}.

Simulated annealing (SA) is a heuristic approach that is widely used for complicated combinatorial optimization problems. It is an iterative procedure that includes two phases: initialization and cool-down. The initial population of candidate solutions can be generated randomly. Then in each iteration of the cool-down phase, a new generation is randomly created from perturbation around the current candidate solutions to generate a new population. This process ends when the termination criteria are satisfied.
In \cite{b22}, SA was used to solve relatively large-scale ridesharing optimization problems.

\section{Proposed method: Biogeography-Based Optimization}\label{Sec:Proposed method}
In this section, we introduce the implementation of our proposed method which leverages an evolutionary approach to search for the globally optimal solution. We found that biography-based optimization (BBO) is an approach that can be effectively applied to the ridesharing problem.

\subsection{Outline of the BBO Algorithm}

Biogeography-based optimization (BBO) \cite{b23} is based on sharing suitable island features. Each island is considered as a possible solution to the problem. Islands, which are encoded as possible solutions for the problem, gradually evolve by sharing features with each other to become more suitable habitats until an acceptably suitable habitat is found, which represents an acceptable solution to the optimization problem. BBO has two important stages, migration and mutation, which are based on the theory of biogeography. Algorithm \ref{algor_bbo_basic} shows a basic implementation of biogeography-based optimization. In the mutation stage, the diversity among habitats is increased, which may result in losing the best solution or finding a better solution. To prevent losing the best solution in each iteration, there are some elites that preserve the current best solution for the next generation.

\begin{algorithm}
\SetAlgoLined
\textbf{Input:} $N$ \;
 Initialize population of candidate solutions $\{x_k\}$ for $k \in [1,N]$\; 
 \While{not termination criterion}{
    \For {each $x_k$} {
        Set emigration probability $\mu_k \propto$ fitness of $x_k$, with $\mu_k \in [0,1]$\;
    }
    \For {each individual $x_k$} {
        Set immigration probability $\lambda_k = 1 - \mu_k$\;
    }
    Sort initial solutions according to cost\;
    Save best solutions as elites\;
    ${z_k} \leftarrow {x_k}$\;
    \For {each individual $z_k$} {
        \For {each solution feature $s$} {
            Use $\lambda_k$ to probabilistically decide whether to immigrate to $z_k$\;
            \If{immigration}{
                Use $\big[ \mu_i \big]_{i=1}^N$ to probabilistically select emigrating individual $x_j$\;
                $z_k(s) \leftarrow x_j(s)$\;
            }
        }
        Probabilistically mutate $\{z_k\}$\;
    }
    $\{x_k\} \leftarrow \{z_k\}$\;
    Replace the worst solutions with elites\; 
 }
\caption{BBO algorithm}
\label{algor_bbo_basic}
\end{algorithm}

The inputs to the algorithm are number of habitats, or population size ($N$); number of elites ($M$); and number of generations ($G_{max}$). In line 1, $N$ solutions (habitats) are randomly initialized. Lines 3 through 6 initialize the emigration rate and the immigration rate for each solution. Lines 10 through 16 use immigration probability to decide on migration for each feature of each solution in the current generation. Finally, in line 18, elites replace the worst solutions in order to keep the previous best solutions from one generation to the next.

\subsection{Adapting BBO for the ridesharing problem}

\begin{figure}
      \centering
      \includegraphics[scale=0.25]{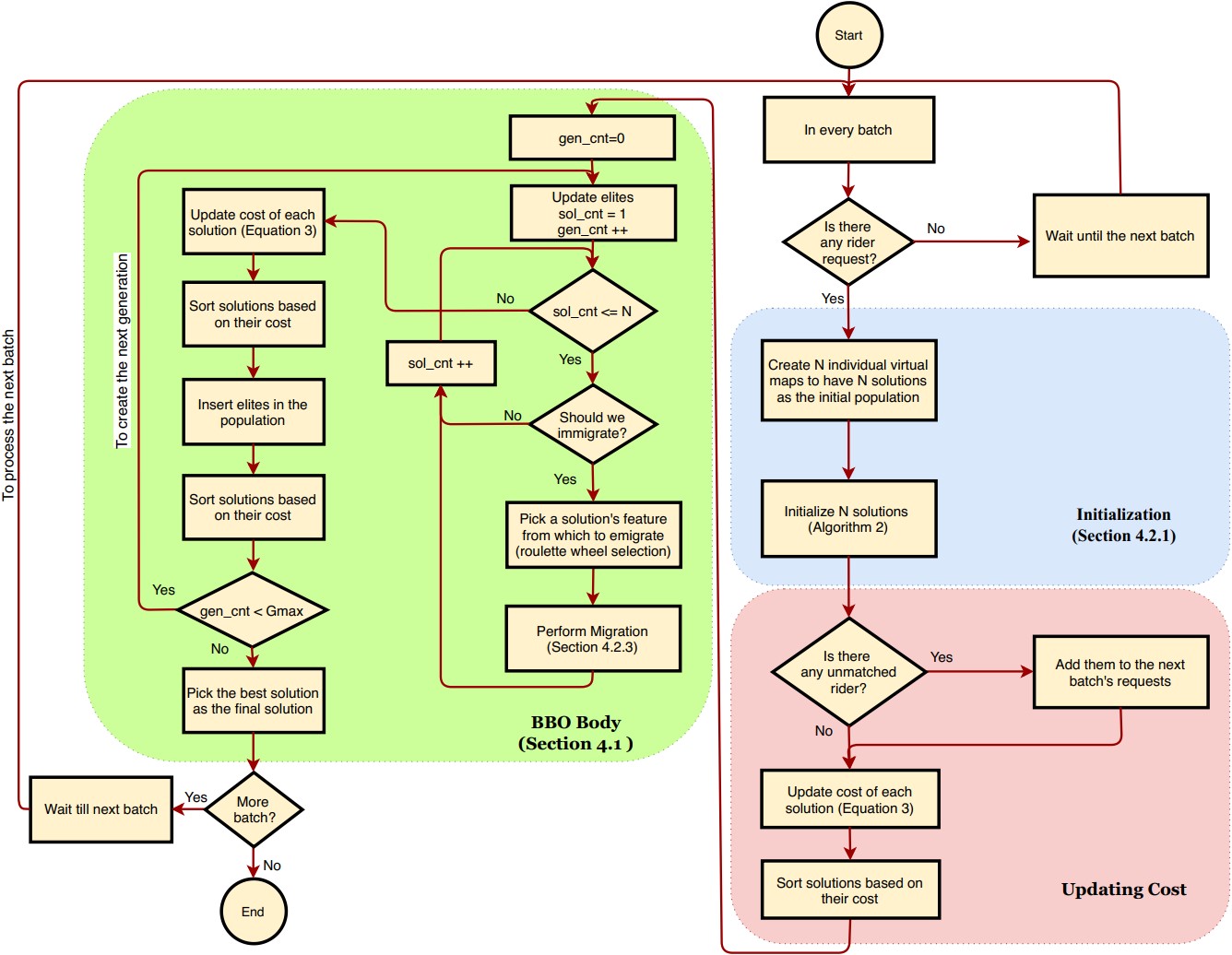}
      \caption{Flowchart of ridesharing optimization using BBO}
      \label{match}
\end{figure}

BBO can be adapted to the ridesharing problem by defining each island as a set of matched trips that include drivers and their assigned riders. The cost of each solution is calculated based on our objective function in Equation \ref{eq:cost}. 
Our method is based on periodic optimization which means we add all new riders to a pool called the rider pool,  and when the period ends, we perform optimization for the rider pool.
We updated the basic BBO algorithm and added some features to it. We customized the initialization and migration steps so that BBO can work with the ridesharing optimization problem.
To handle rider requests, we gather requests which appear within a certain time slot called a batch, and then try to solve each batch like a static ridesharing problem. In the first step of each batch, we generate the initial population, and in the second step we update the cost of each solution based on the cost function. We iteratively evolve the initial population to generate better solutions with lower cost in each successive population. Figure \ref{match} shows the outline of our proposed method. 

\subsubsection{Initialization}

The first step of the BBO algorithm is initialization, which means $N$ candidate solutions need to be generated. In our simulation, which will be described in Section \ref{Sec:Experiments and results}, we put drivers and riders on the map based on their origin positions as they become available. The simulator updates the positions of the drivers on the map based on their schedules. In the next section, we show how each solution is stored and represented in our framework.

Algorithm \ref{algorinit} shows our proposed initialization function. The first step in this function is making the initial population. In order to have some good initial solutions in the initial generation, we design the initialization function so that it generates some solutions randomly and some solutions more intelligently, such as with a greedy algorithm. $H_{ratio}$ defines what proportion of the initial population will be generated randomly. If $H_{ratio} = 0$, the entire initial population will be generated with the greedy algorithm. 

One  of the important characteristics of a good initial population is that all solutions should be unique and widely scattered within the problem's domain. In order to achieve a diversity of initial solutions, the initialization algorithm randomly selects a rider from each batch as a starting point so that the greedy algorithm generates a different solution each time it executes.  

In Algorithm \ref{algorinit}, MSP is the minimum shortest path that a driver needs to travel. The $MSP$ of a vehicle's schedule in a weighted graph is a route with minimum total weight. $MSP$ in a weighted graph of n vertices can be found in time $O(n^3)$.
To have a better view on the proposed method, Figure \ref{match} is provided which shows how our dynamic ridesharing framework works.

\begin{algorithm}
\SetAlgoLined 
    Create $N$ virtual maps from the original map \;
    Initialize population of candidate solutions ${x_k}$ for $k \in [1,N]$\; 
    Apply each initial solution to a virtual grid\;
    \For{$ k $ = 1 to $N$ }{
        \For {each $ r \in R $}{
            $cand\_list \leftarrow $ feasible drivers $\subset D$\;
            \eIf{$ k < N H_{ratio} $ }{
                //***Greedy assignment***\;
                $d^* \leftarrow $ driver from $cand\_list$ with minimum $d_{ov}$ for rider $r$\; 
                \eIf{Constraints of $d^*$ and $r$ are satisfied }{
                    Assign $d^*$ to rider $r$ and update schedule for $d^*$ \;
                    Update $d^*_{rte} \leftarrow \mbox{MSP}(schedule) $\;
                }
                {
                    Remove $d^*$ from $cand\_list$ and goto Greedy assignment 
                }
            }
            {
                //***Random assignment***\;
                $d^* \leftarrow$ random driver from $cand\_list$\;
                \eIf{Constraints of $d^*$ and $r$ are satisfied }{
                    Assign $d^*$ to rider $r$ and update schedule for $d^*$ \;
                    Update $d^*_{rte} \leftarrow MSP(schedule) $\;
                }
                {
                    Remove $d^*$ from $cand\_list$ and goto Random assignment 
                }
            }
        }
    }
    \caption{Generate Initial Population}
    \label{algorinit}
\end{algorithm}

\subsubsection{Problem representation}
We need $N$ initial solutions with the same map, same drivers, and same riders. Due to simulator limitations we cannot define multiple solutions for a single map. We need $N$ maps and $N$ solutions that can be processed simultaneously. So we define the concept of virtual maps, which are just duplicates of the initial map with objects on it such as drivers and riders. In this way, we can generate a solution for $N$ virtual maps and then update each virtual map's objects with the corresponding solution. Thus, we copy the initial map to $N$ virtual maps in order to perform optimization, and once optimization is done we pick the best map (best solution) and set it as the initial map for the next batch. 

\begin{figure}
    \centering
    \includegraphics[scale=0.21]{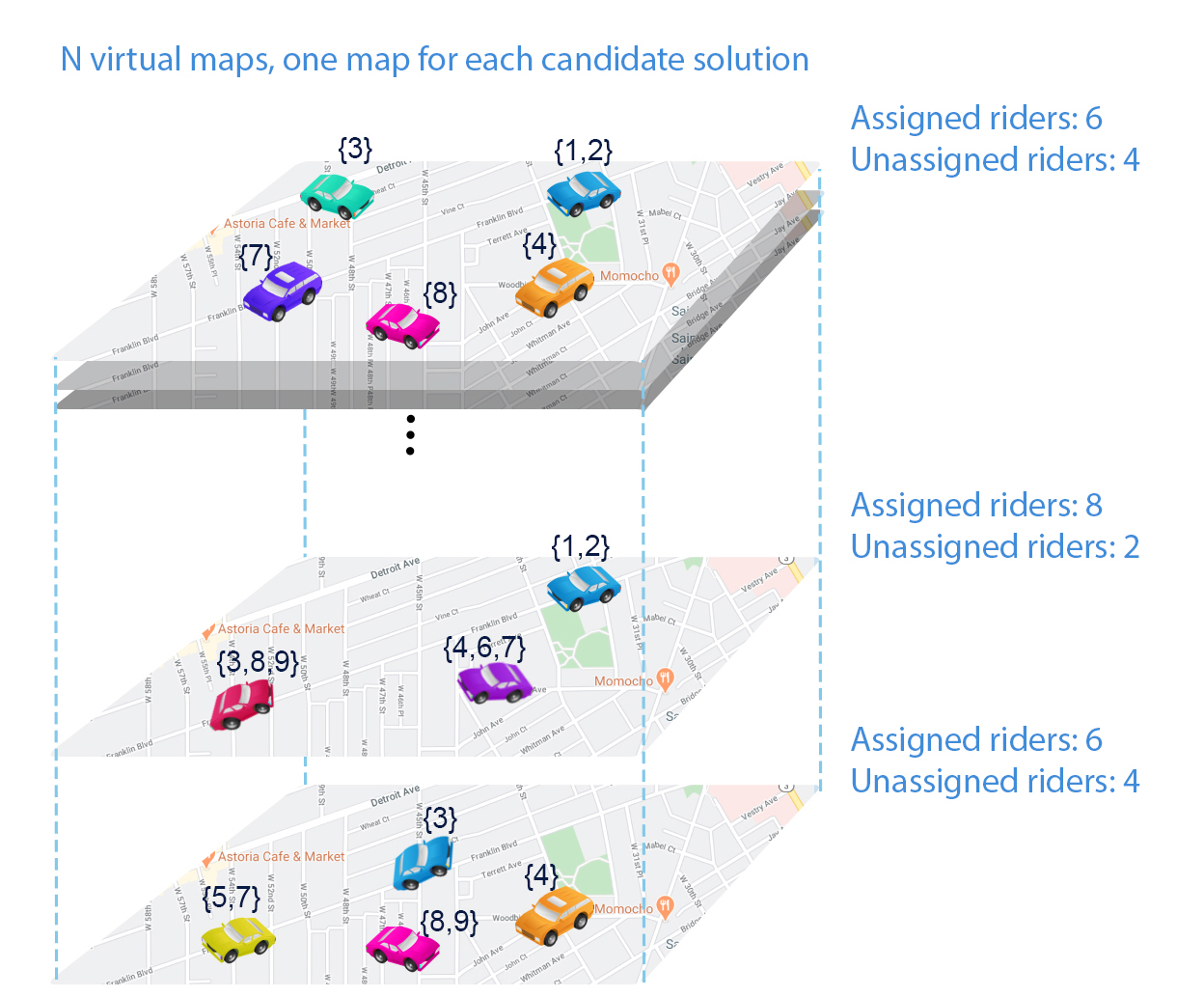}
    \caption{Concept of virtual maps. This figure shows $N$ candidate solutions during a given batch. Each candidate solution represents a different possible driver/rider matching strategy.}
    \label{Vmap}
\end{figure}

Figure \ref{Vmap} shows an example using virtual maps. There are $N$ virtual maps corresponding to $N$ solutions, of which only three are shown in the figure. There are many drivers on each map, but the figure only shows those drivers that are assigned to one or more riders in the current batch. It is necessary to generate multiple unique initial solutions to cover as much as possible of the solution space. Therefor, our initialization method has to generate unique solutions as shown in Figure \ref{Vmap}. As seen in the figure, there are 10 request in the rider pool (the total number of assigned riders and unassigned riders in each candidate solution). The solution corresponding to the virtual map in the middle can be considered as a good solution since it includes more matches for the riders than the other solutions. The numbers shown above each vehicle represents the assigned riders in that vehicle.

\subsubsection{Migration approach}
\label{sec-Migration}
\begin{figure}
      \centering
      \includegraphics[scale=0.33]{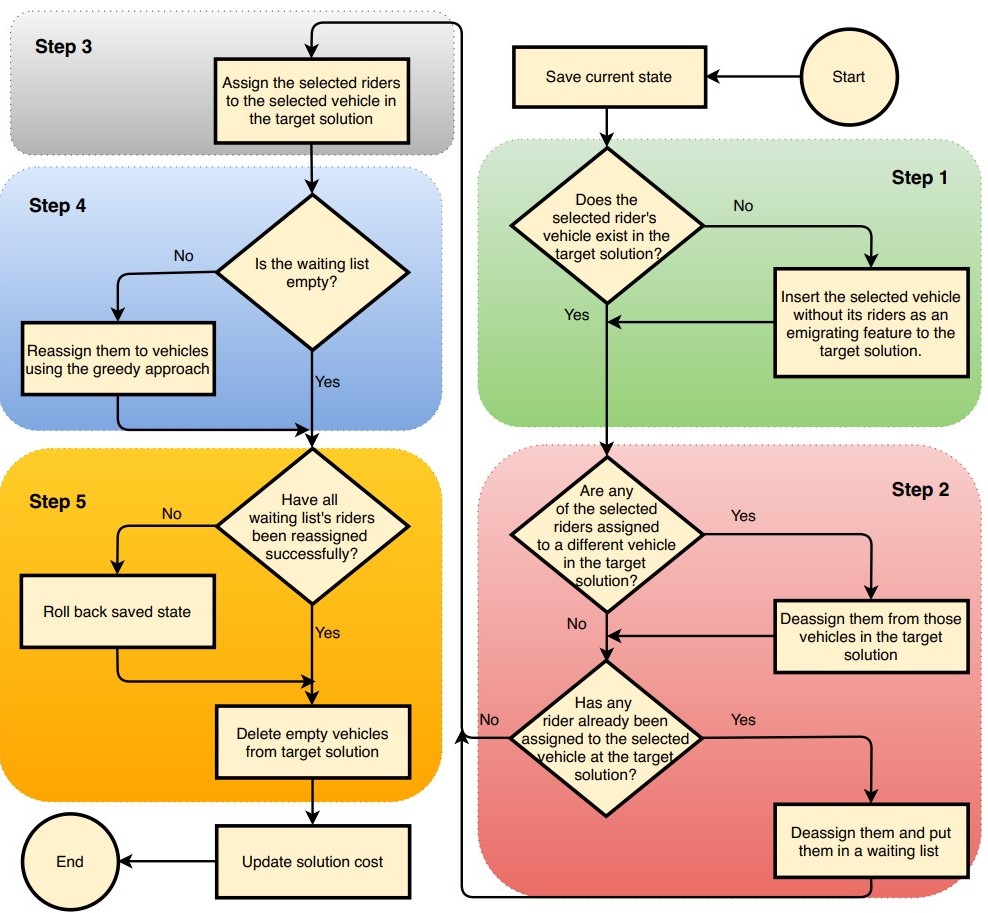}
      \caption{Flowchart of migration function}
      \label{migration}
\end{figure}

In this section, the process of how to migrate features between two islands in BBO is described. So, we explain what are the features and what are the islands in the ridesharing optimization problem. In our proposed migration approach we define all matched drivers on a virtual map as determined by the initialization step for the island. Then, we define each vehicle in the set of matched drivers as a feature of that island. Figure \ref{migration} depicts the details of this approach, which consists of the following steps.
\begin{enumerate}
\item Insert the selected rider's vehicle (from the emigrating solution) in the target solution (into the immigrating solution) if it does not yet exist.
\item If the selected vehicle already exists in the target solution, remove its current riders and put them in a waiting list.
\item Assigned the selected riders to the selected vehicle.
\item Reassign all riders that are in the waiting list to new drivers using a greedy approach.
\item If any rider from the waiting list remains unassigned, roll back (undo) this migration sequence.
\end{enumerate}

\section{Experiments and Results}\label{Sec:Experiments and results}
In this section, we introduce the dataset that is used in this research and then we show experimental results obtained by our proposed method. Finally, we compare our results with other algorithms, including one state-of-the-art method.

\subsection{Dataset}
We conduct all experimental evaluations on a real Beijing road network \cite{b15}. This dataset was designed for taxi-ridesharing simulations. Since our research is focused on peer-to-peer ridesharing problems, in which each driver has its own specific destination and time constraints, the dataset is not directly applicable and must be changed.
We modified this dataset so that each vehicle has a destination, which is a random location on the map, and a certain deadline to reach the destination. The deadline is set so that each driver has twice the minimum required time (MSP) to reach his destination.
Our dataset consist of two parts as follows.
\begin{enumerate}
    \item{Road network:} There are 351,290 vertices and 743,822 edges in a total area of 17,158 km$^2$ in the Beijing road network dataset. Each vertex contains latitude and longitude. For each directly connected edge, the travel distance is given. The shortest path distance between any two vertices can be obtained by searching an undirected graph of the road network.
    
    \item{Problem instance:} Our problem instance is flexible and can be defined as a small problem or a large problem in order to evaluate the scalability and performance of a given algorithm in different circumstances. The number of drivers, number of riders, and their spatial distributions are important factors that we can change based on the purpose of the experiment. Although the number of drivers and riders can be manipulated, some statistical facts \ie the distribution of the rider appearances, are based on Chinese transportation data. In the dataset file, each record is a driver or a rider, which can be determined based on the $ load $ feature. $ load $ is a negative number for drivers and a positive number for riders. Each driver has the properties id ($d_{id}$), origin ($d_o$), early time ($d_e$), late time ($r_l$), and capacity ($d_{cap}$). Each rider has the properties id ($r_{id}$), origin ($r_o$), destination ($r_d$), early time ($r_e$), and late time ($r_l$). 
    \begin{enumerate}
        \item{Driver:} The initial locations of the drivers are sampled from the real trip dataset until the desired number of drivers is obtained. Each sampled trip origin is used as the initial location of the driver. Each vehicle has three seats in addition to the driver's seat, and is initialized as a vertex on the road network. When a driver is initialized, it has no riders and there is only one location in its schedule, which is the driver's origin. The driver's destination is selected randomly and is added to its schedule later. 
        \item{Rider:} Riders are extracted from the real dataset from all taxi trips in the sample period from 6:00–6:30 PM on an arbitrary day in the Beijing network. There are 17,467 trips, about 9.70 requests per second. The request rate can be changed to simulate customer appearance frequencies of 0.5x, 2.0x, and 4.0x.
        \item{Time Window:} Each rider has own early time constraint and late time constraint.
    \end{enumerate}

\end{enumerate}

\subsection{Simulator}

We used the Cargo ridesharing simulator to implement our ridesharing optimization algorithms \cite{b15}. Cargo is a simulator that provides fundamentals of a ridesharing system, \ie vehicle motion, spatial indices, statistics, and event logs.

\begin{figure}
    \centering
    \includegraphics[scale=0.45]{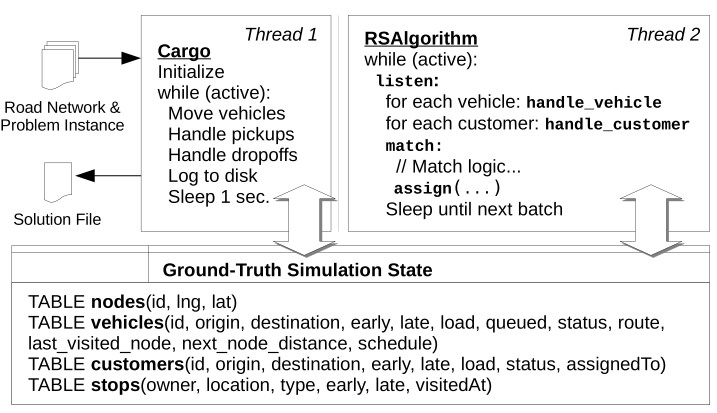}
    \caption{Cargo architecture \cite{b15}}
    \label{Cargo_arch}
\end{figure}

One of the important features of Cargo is that it accepts as an input a matching algorithm, which determines the driver / rider matching policies and objectives. Thus, we implemented our algorithms and attached them to Cargo. Figure \ref{Cargo_arch} shows the details of the Cargo architecture. Cargo and RSAlgorithm are two important components of the simulator. Cargo is the core component of the simulator and maintains the problem instance and the road network. It also is responsible for creating statistics, logs, and reports. RSAlgorithm provides the outline of the ridesharing algorithm. It provides events that can be used in any user-defined algorithm. 

To compute the shortest path between two locations, this simulator uses the G-tree algorithm, which is a state-of-the-art method for quickly computing the minimum path. The G-tree file index of the Beijing road network is pre-built and is available as part of Cargo.

\subsection{Results}
The dataset provided in Cargo is a real-world dataset, so in order to ease code development and debugging, we made a smaller problem from the existing dataset. In the new dataset, there are 668 riders and 200 drivers that appear throughout a time period of 30 minutes. 
We set the time of each batch to 30 seconds and we optimize matching during each batch, so 60 batches comprise 30 minutes of simulation time. Riders enter the ridesharing system with a uniform distribution so each of the first sixty 30-second batches include 11 new riders, and the last batch includes 8 riders to complete the total of 668 riders. There are 22 drivers who appear during the first batch, three drivers who appear in each of the following 59 batches, and 1 driver who appears in the last batch which results in a total of 200 drivers.
To demonstrate the impact of different parameters on the performance of our BBO method, we test our proposed algorithm for four different sets of BBO parameters as shown in Table \ref{tab:cases}. To have a fair comparison, each case uses the same initial population. Table \ref{tab:cases_results} shows the performance of each of our four test cases. In Table \ref{tab:cases_results}, \textit{base driver distance} is the sum of all the drivers' own trip distances; \textit{base rider distance} is the sum of all the riders' own trip distances, i.e., the distance from each rider's origin to its destination; and \textit{matched trip distance} is the total distance of all trips after matching riders with drivers. \textit{Base driver distance} and \textit{base rider distance} are the same for all cases since these numbers depend only on the input dataset. Table \ref{tab:cases_results} also shows $d_{ov}$ and $M_R$, which are important parts of the objective function; and \textit{cost}, which is computed based on the objective function definition in Equation \ref{eq:cost}.

\begin{table}[]
\centering
\caption{BBO parameters for four different sets of parameters (see Algorithm \ref{algor_bbo_basic} for parameters)}
\label{tab:cases}
\begin{tabular}{llllll}
\hline
                        & Symbol & Case 1 & Case 2 & Case 3 & Case 4 \\ \hline
Generation limit       & $G_{max}$   & 10     & 10     & 10     & 10     \\
Population size         & $N$      & 20     & 20     & 20     & 20     \\
Hybrid  init. rate & $H_{ratio}$  & 0.85   & 0.85   & 1      & 1      \\
Number of elites            & $E$      & 1      & 1      & 1      & 1      \\
Roll back                   & $RB$     & On     & Off    & On     & Off    \\
Objective weight            & $\alpha$    & 0.5    & 0.5    & 0.5    & 0.5    \\
Cache size & - & 2\textasciicircum{}30 & 2\textasciicircum{}30 & 2\textasciicircum{}30 & 2\textasciicircum{}30 \\ \hline
\end{tabular}
\end{table}

\begin{table}[]
\centering
\caption{Results of the proposed BBO algorithm. The best result in each row is shown in bold font.}
\label{tab:cases_results}
\begin{tabular}{lllll}
\hline
                                            & Case 1      & Case 2      & Case 3      & Case 4      \\ \hline
Base driver dist.                   & 2,168,819   & 2,168,819   & 2,168,819   & 2,168,819   \\
Base rider dist.                    & 3,532,517   & 3,532,517   & 3,532,517   & 3,532,517   \\
Matched trip dist.                  & 3,483,903   & 3,496,608   & 3,454,594   & {\bf 3,451,805}   \\
$d_{ov}$                                       & 1,315,084   & 1,327,789   & 1,285,775   & {\bf 1,282,986}    \\
$M_R$                                       & 0.479       & {\bf 0.482}       & 0.481       & 0.476       \\
Cost                                        & 0.8932      & 0.8938      & {\bf 0.8834}      & 0.8871      \\ \hline
\end{tabular}
\end{table}

\begin{figure*}%
    \centering
  
    \subfloat[\centering Matching rate  \label{Matchingrate} ]{{\includegraphics[width=5cm]{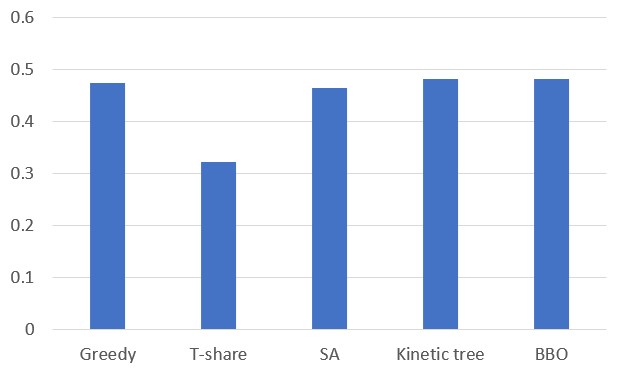} }}%
    \qquad
    \subfloat[\centering Distance overhead \label{dov}]{{\includegraphics[width=5cm]{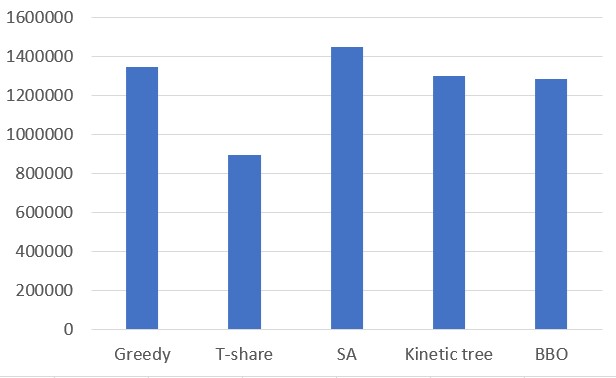} }}%
    \qquad
    \subfloat[\centering cost \label{Cost}]{{\includegraphics[width=5cm]{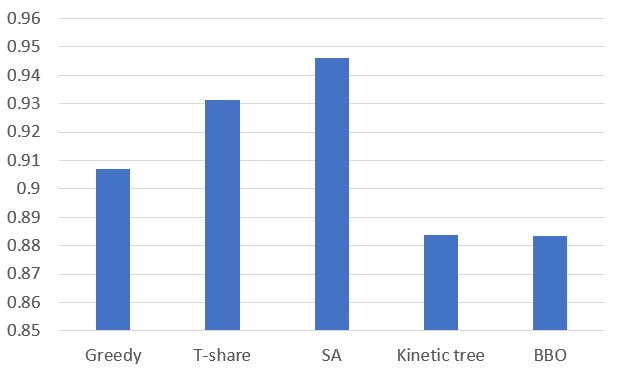} }}%
    
    \caption{Algorithms' performance comparison }%
    \label{results_comp}%
\end{figure*}

Table \ref{tab:cases_results} shows approximately the same performance for all cases; however, the minimum cost belongs to case 3, which uses all-random initial solutions and takes advantage of the \textit{roll back} mechanism (see Section~\ref{sec-Migration}). Although case 3 has neither the highest $M_R$ nor the lowest $d_{ov}$, it has the minimum cost when both objectives are combined based on the multi-objective cost function. 

The other thing that can be inferred from Table \ref{tab:cases_results} is that using a 100\% random initial population results in a lower cost than using more sophisticated initial solutions, e.g., greedy-based solutions, in the initial population. Table \ref{tab:cases_results} also shows there is essentially no advantage or disadvantage if we use \textit{roll back}.

To evaluate our proposed BBO method with other research, we implemented four other approaches in the simulator: greedy algorithm (GR), simulated annealing (SA), T-share (TS) \cite{b11} and kinetic tree (KT) \cite{b18}. In order to have a fair comparison among these algorithms and our BBO method, we use the same dataset and the same simulation parameters for all methods.
We compare all algorithms in terms of $d_{ov}$, $M_R$, and cost in Figures \ref{Matchingrate}--\ref{Cost}.



Figure \ref{results_comp} \subref{Matchingrate} shows that the highest $M_R$ corresponds to BBO and KT with 0.481 and 0.482 respectively, with the Greedy and SA algorithms very close behind, and the lowest matching rate belongs to TS with only 0.322. Figure \ref{results_comp} \subref{dov} shows that TS has the best distance overhead and SA has the worst. The reason that TS has the minimum $d_{ov}$ is that it matches fewer riders with drivers, so drivers have fewer matched riders on-board, which results in less distance overhead. Finally, Figure \ref{results_comp} \subref{Cost} shows that the best cost belongs to BBO and KT with 0.8834 and 0.8837 respectively, while SA has the worst cost.


\section{Conclusion and Future Work}
We proposed an algorithm which leverages biogeography-based optimization to find a near-optimal solution for the online ridesharing problem while imposing minimum delay for driver / rider matching. We evaluate our proposed algorithm using a Cargo simulator on a Beijing road network. The proposed algorithm performs competitively with the well-known simulated annealing, and T-share,  and kinetic tree algorithms.

To improve the performance of our proposed algorithm, we aim to implement BBO in a distributed fashion. In other words, we will partition the map into smaller domains and then optimize each domain individually, but with each domain being aware of its neighbors solutions.

\end{document}